\documentclass[12pt]{amsart}
\usepackage{amsmath,amscd}
\usepackage{amssymb}
\usepackage[arrow,matrix,curve]{xy}

\title{Actors in categories of interest}
\author[J. M. Casas, T. Datuashvili and M. Ladra]{J. M . Casas$^1$, T. Datuashvili$^2$ and M. Ladra$^3$}
 \thanks{The first and third authors were
supported by the MEC (Spain), Grant MTM 2006-15338-C02 (European
Feder support included). The second author is grateful to Santiago
and Vigo Universities for financial supports.}

\newcommand{\vf}{\varphi}

\newcommand{\tht}{\theta}
\newcommand{\pa}{\partial}
\newcommand{\bC}{\mathbb{C}}

\newcommand{\bB}{\mathbb{B}}
\newcommand{\cB}{\mathfrak{B}}
\newcommand{\bE}{\mathbb{E}}
\newcommand{\bJ}{\mathbb{J}}
\newcommand{\bL}{\mathbb{L}}
\newcommand{\bG}{\mathbb{G}}
\newcommand{\Om}{\Omega}
\newcommand{\om}{\omega}
\newcommand{\aaa}{{\bf{\mathfrak{a}}}}
\newcommand{\bbb}{{\bf b}}

\newcommand{\lra}{\longrightarrow}
\newcommand{\ol}{\overline}
\newcommand{\wt}{\widetilde}

\newcommand{\os}{\overset}
\newtheorem{theo}{Theorem}[section]
\newtheorem{Lem}[theo]{Lemma}
\newtheorem{Prop}[theo]{Proposition}

\newtheorem{defi}[theo]{Definition}
\newtheorem{axiom}{Axiom}
\newtheorem{cond}{Condition}
\newtheorem{remark}{Remark}
\newtheorem{example}{Example}
\numberwithin{equation}{section}

\begin{document}
\maketitle

\noindent{ $^1$ \footnotesize \it Dpto. Matem\'atica Aplicada
I, Univ. de Vigo, 36005 Pontevedra, Spain, e-mail: jmcasas@uvigo.es}\\
{$^2$ \footnotesize \it A. Razmadze Math. Inst.,
Alexidze str.1, 0193 Tbilisi,  Georgia, e-mail: tamar@rmi.acnet.ge}\\
{$^3$ \footnotesize \it Departamento de \'Algebra,  Universidad de
Santiago de Compostela, Santiago de Compostela 15782, Spain, e-mail:
ladra@usc.es}

\begin{abstract}
For an object $A$ of a category of interest $\mathbb{C}$ we
construct the group with operations ${\mathfrak{B}}(A)$ and the
semidirect product ${\mathfrak{B}}(A)\ltimes A$ and prove that there
exists an actor of $A$ in $\mathbb{C}$ if and only if
${\mathfrak{B}}(A)\ltimes A \in \mathbb{C}$. The cases of groups,
Lie, Leibniz, associative, commutative associative, alternative
algebras, crossed and precrossed modules are considered. In
particular, we give examples of categories of interest, where always
exist actors. The paper contains some results for the case $\Omega_2
= \{+,*,*^{\circ}\}$.
\end{abstract}

\section{Introduction}
 The paper is dedicated to the question of
the existence and construction of actors for the objects in
categories of interest (see below Section 2 for the definitions
and examples). This kind of categories were introduced by G.
Orzech \cite{Or}. Actions in
  algebraic categories were studied by G. Hochschild \cite{Ho}, S. Mac Lane \cite{Mc}, A. S.-T. Lue \cite{Lu}, K. Norrie \cite{No},
J.-L. Loday \cite{Lo},
  R. Lavendhomme and Th. Lucas \cite{LL} and others. The authors were looking for the
  analogs of automorphisms of groups in associative algebras, rings, Lie
  algebras, crossed modules and Leibniz algebras. We see different approaches to this question.
  Lue and Norrie (on the base of the results of Lue \cite{L1} and Whitehead \cite{Wh}),
  to any object associate a certain type of object, the construction
  in the corresponding category, called actor of this object \cite{No}, with
  special properties, analogous to group automorphisms, under which
  is meant that the actor fits into a certain commutative diagram (see
  Section 2, diagram (2.7)). In \cite{LL} Lavendhomme and Lucas introduce
  the notion of a $\Gamma$-algebra of derivations for an algebra $A$,
  which is the terminal object in the category of crossed modules
  under $A$. Recently F. Borceux, G. Janelidze and G. M. Kelly
  \cite{BJK}
  proposed a categorical approach to this question. They study internal
  object actions  and introduce the notion of a
  representable action, which in the case of a category of interest
  is equivalent to the definition of an actor given in this paper (see Section 3).

Let $\bC$ be a category of interest with a set of operations
$\Om=\Om_0\cup\Om_1\cup\Om_2$ and a set of identities $\bE$. We
define a general category of groups with operations $\bC_G$ with the
same set of operations and a set of identities
$\bE_G\hookrightarrow\bE$ for which $\bC\hookrightarrow\bC_G$. We
introduce the notions of an actor and of a general actor object for
the objects of $\bC$. For any object $A\in\bC$ we give a
construction of the universal algebra $\cB(A)$ with the operations
from $\Om$. We show that in general $\cB(A)$ is an object of
$\bC_G$. For any $A\in\bC$ we define an action of $\cB(A)$ on $A$,
which is a $\cB(A)$-structure on $A$ in $\bC_G$ (i.e. the derived
action appropriate to $\bC_G$, see Section 2 below for the
definitions). In a well-known way we define the universal algebra
$\cB(A)\ltimes A$ which is an object of $\bC_G$. We define the
homomorphism $A\lra \cB(A)$ in $\bC_G$, which turned out to be a
crossed module in $\bC_G$. We show that the general actor object
always exists and $\cB(A)$=GActor$(A)$ (Theorem 3.7). The main
theorem states that an object $A$ from $\bC$ has an actor in $\bC$
if and only if $\cB(A)\ltimes A$ is an object in $\bC$ and in this
case $A\lra \cB(A)$ is an actor of $A$ in $\bC$ (Theorem 3.6). From
the results of \cite[Theorem 6.3]{BJK} and from Theorem 3.6 of this
paper we conclude that a category of interest $\mathbb{C}$ has
representable object actions in the sense of \cite{BJK} if and only
if $\cB(A) \ltimes A \in \mathbb{C}$ for any $A \in \mathbb{C}$, and
if it is the case the corresponding representing objects are
$\cB(A), A \in \mathbb{C}$.

 We consider separately the case $\Om_2=\{+,*,*^\circ\}$. In the case of groups
$(\Om_2=\{+\})$ we obtain that $\cB(A)\approx \text{Aut}(A)$,
$A\in\bG$r. In the case of Lie algebras $(\Om_2=\{+,[\;,\;]\})$ for
$A\in \bL$ie we obtain $\cB(A)\approx\text{Der}(A)$. In the case of
Leibniz algebras we have $\cB(A)\in\bL$eib, for any $A\in\bL$eib;
$\cB(A)$ has a derived set of actions on $A$ if and only if for any
$B,C\in\bL$eib which have a derived action on $A$ we have
$[c,[a,b]]=-[c,[b,a]]$, for any $a\in A$, $b\in B$, $c\in C$, (which
we call Condition 1 and it is equivalent to the existence of an
Actor$(A)$). In this case $\cB(A)=\text{Actor}(A)$ (Theorem 4.5). We
give examples of such Leibniz algebras. In particular Leibniz
algebras $A$ with Ann$(A)=(0)$, where Ann$(A)$ denotes the annulator
of $A$, and  perfect Leibniz algebras (i.e. $A = [A, A]$) satisfy
Condition 1. We have an analogous picture for associative algebras.
In this case $\cB(A)$ is always an associative algebra, but the
action of $\cB(A)$ on $A$ defined by us is not a derived action on
$A$. Here we introduce Condition 2: for any $B$ and $C \in
\mathbb{A}$ss, which have a derived action on $A$ we have $c *\ (a*\
b)=(c*\ a)*\ b$, for any $a\in A$, $b\in B$, $c\in C$, where $*$
denotes the action. The action of $\cB(A)$ on $A$ is a derived
action if and only if $A$ satisfies Condition 2 and it is equivalent
to the existence of an Actor$(A)$. In this case
$\cB(A)=\text{Actor}(A)$ (Proposition 4.6). Associative algebras
with conditions Ann$(A)=(0)$ or with $A^2=A$ satisfy Condition 2.
These kind of associative algebras are considered in \cite{LL,Mc}).
The cases of modules over some ring, commutative associative
algebras, alternative algebras, crossed modules and precrossed
modules in the category of groups are discussed. Note, that the
construction and the results given in our work enabled us to prove
the existence of an actor in the category of precrossed modules. We
consider the case where $\Om_2=\{+,* ,*^\circ\}$. The necessary and
sufficient conditions for the existence of an actor is determined in
the case where $\bE$ contains only all identities from $\bE_G$ and
Axiom 1 and Axiom 2, and Axiom 2 has not consequence identities
(Theorem 4.11). We consider separately the algebras of the same type
as in Theorem 4.11 with additional commutativity or
anticommutativity condition. We obtain the necessary and sufficient
conditions for the existence of an actor in the corresponding
category (Theorem 4.12) and give examples of such categories. The
paper contains a comment to the formulation of Proposition 1.1 of
\cite{Da}, which we apply in this paper. The results obtained here
enabled us to find the construction of an actor of a precrossed
module, which will be the subject of the forthcoming paper. The
results of the paper are included in the Doctoral dissertation of
the second author \cite{Da2} Chapter V.

In   Section 2 we present the main definitions and results which are
used in what follows. In Section 3 we give the main construction of
the object $\cB(A)$ and the corresponding results. In Section 4 we
consider the case of groups, Lie algebras, associative algebras and
Leibniz algebras. For the special types of objects in $\mathbb{A}$ss
and $\bL$eib it is proved that $\cB(A)\approx\text{Bim}(A)$,
$\cB(A)\approx\text{Bider}(A)$ respectively, (Propositions 4.7,
4.8), where Bim$(A)$ denotes the associative algebra of
bimultipliers defined by G. Hochschild and by S. Mac Lane for rings
(called bimultiplications in \cite{Mc} and multiplications in \cite
{Ho} from where the notion comes) and Bider$(A)$ denotes the Leibniz
algebra of biderivations of $A$ defined in Section 2, which is
isomorphic for these special types of Leibniz algebras to the
biderivation algebra defined by J.-L. Loday in \cite{Lo}. At the end
of the section we summarize for the special type of categories of
interest the results of Sections 3 and 4 and obtain the necessary
and sufficient conditions for the existence of an actor.

\section{Preliminary definitions and results}

This section contains well-known definitions and results which
will be used in what follows.

Let $\bC$ be a category of groups with a set of operations $\Om$
and with a set of identities $\bE$, such that $\bE$ includes the
group laws and the following conditions hold. If $\Om_i$ is the
set of $i$-ary operations in $\Om$, then:

(a) $\Om=\Om_0\cup\Om_1\cup\Om_2$;

(b) The group operations (written additively : $0,-,+$) are
elements of $\Om_0$, $\Om_1$ and $\Om_2$ respectively. Let
$\Om_2'=\Om_2\setminus\{+\}$, $\Om_1'=\Om_1\setminus\{-\}$ and
assume that if $*\in\Om_2$, then $\Om_2'$ contains $*^\circ$
defined by $x*^{\circ}y=y*x$. Assume further that $\Om_0=\{0\}$;

(c) for each $*\in\Om_2'$, $\bE$ includes the identity
$x*(y+z)=x*y+x*z$;

(d) for each $\om\in\Om_1'$ and $*\in\Om_2'$, $\bE$ includes the
identities $\om(x+y)=\om(x)+\om(y)$ and $\om(x)*y=\om(x*y)$.

We formulate more axioms on $\bC$ (Axiom (7) and Axiom (8) of
\cite{Or}).

If $C$ is an object of $\bC$ and $x_1,x_2,x_3\in C$:


\begin{axiom}
 $x_1+(x_2*x_3)=(x_2*x_3)+x_1$, for each
$*\in\Om_2'$.
\end{axiom}

\begin{axiom} For each ordered pair
$(*,\ol{*})\in\Om_2'\times\Om_2'$ there is a word $W$ such that
\begin{gather*}
(x_1*x_2)\ol{*}x_3=W(x_1(x_2x_3),x_1(x_3x_2),(x_2x_3)x_1,\\
(x_3x_2)x_1,x_2(x_1x_3),x_2(x_3x_1),(x_1x_3)x_2,(x_3x_1)x_2),
\end{gather*}
where each juxtaposition represents an operation in $\Om_2'$.
\end{axiom}

 We will write the right side of Axiom 2 as
$W(x_1,x_2;x_3;*,\ol{*})$. A category of groups with operations
satisfying  Axiom 1 and Axiom 2 is called a category of interest
by Orzech \cite{Or} (see also \cite{Po}).

Note that from the equalities
$(x+y)*(z+t)=x*z+x*t+y*z+y*t=x*z+y*z+x*t+y*t$ follows that
$x*t+y*z=y*z+x*t$, for $*\in\Om_2'$, $x,y,z,t\in C$, $c\in\bC$.

Denote by $\bE_G$ the subset of identities of $\bE$ which includes
the group laws and the identities (c) and (d). We denote by $\bC_G$
the corresponding category of groups with operations. Thus we have
$\bE_G \hookrightarrow \bE$, $\bC=(\Om,\bE)$, $\bC_G=(\Om,\bE_G)$
and there is a full inclusion  functor $\bC\hookrightarrow\bC_G$.

In the case of associative algebras with multiplication
represented by  $*$, we have $\Om_2'=\{*,*^\circ\}$. For Lie
algebras take $\Om_2'=([\;,\;],[\;,\;]^\circ)$ (where
$[a,b]^\circ=[b,a]=-[a,b]$). For Leibniz algebras (see  the
definition below), take $\Om_2'=([\;,\;],[\;,\;]^\circ)$, (here
$[a,b]^\circ=[b,a]$). It is easy to see that all these algebras
are categories of interest. In the example of groups
$\Om_2'=\varnothing$. As it is mentioned in \cite{Or} Jordan
algebras do not satisfy Axiom 2.


\begin{defi} \cite{Or} Let $C\in\bC$. A subobject of $C$ is called an ideal if it is the kernel of some morphism.

\end{defi}

\begin{theo}\cite{Or} Let $A$ be a subobject of $B$ in $\bC$. Then $A$ is an ideal of $B$ if and only if the following conditions hold;

{\rm (i)} $A$ is a normal subgroup of $B$;

{\rm (ii)} For any $a\in A$, $b\in B$ and $*\in\Om_2'$, we have
$a*b\in A$.
\end{theo}

\begin{defi}\cite{Or}
Let $A$, $B\in\bC$. An extension of $B$ by $A$ is a sequence
\begin{equation}
\xymatrix{0\ar[r]&A\ar[r]^-{i}&E\ar[r]^-{p}&B\ar[r]&0}
\end{equation}
in which $p$ is surjective  and  $i$ is the kernel of $p$. We say
that an extension is split  if there is a morphism $s:B\lra E$
such that $ps=1_B$.
\end{defi}

\begin{defi}\cite{Or}
A split extension of $B$ by $A$ is called a $B$-structure on $A$.
\end{defi}

As in \cite{Or}, for $A,B\in\bC$ we will say we have ``a set of
actions of $B$ on $A$'', whenever there is a set of maps
$f_*:B\times A\lra A$, for each $*\in\Om_2$.

A $B$-structure induces a set of  actions of $B$ on $A$
corresponding to the operations in $\bC$. If (2.1) is a split
extension, then for $b\in B$, $a\in A$ and $*\in\Om_2{}'$ we have
\begin{align}
b\cdot a &=s(b)+a-s(b),\\
b*a&=s(b)*a.
\end{align}
(2.2) and (2.3) are called derived  actions of $B$ on $A$ in
\cite{Or} and split derived actions in \cite{Da}, since we
considered there actions derived from non split extensions too
when $A$ is a singular object.

Given a set of actions of $B$ on $A$ (one for each operation in
$\Om_2$), let $B\ltimes A$ be a universal algebra whose underlying
set is $B\times A$ and whose operations are
\begin{align*}
(b',a')+(b,a)&=(b'+b,a'+b'\cdot a),\\
(b',a')*(b,a)&=(b'*b,a'*a+a'*b+b'*a).
\end{align*}

\begin{theo}\cite{Or}
A set of actions of $B$ on $A$ is a set of derived actions if and
only if $B\ltimes A$ is an object of $\mathbb{C}$.
\end{theo}

Together with the description of the set of derived actions given
in the Theorem above, we will need the identities which satisfies
a set of derived actions in case $A,B\in\bC_G$ and which guarantee
that the set of actions is a set of derived actions in $\bC_G$.

\begin{Prop}\cite{Da}
A set of actions in $\bC_G$ is a set of derived actions if and
only if  it satisfies the following conditions:

$1$. $0\cdot a=a$,

$2$. $b\cdot(a_1+a_2)=b\cdot a_1+b\cdot a_2$,

$3$. $(b_1+b_2)\cdot a=b_1\cdot(b_2\cdot a)$,

$4$. $b*(a_1+a_2)=b* a_1+b* a_2$,

$5$. $(b_1+b_2)* a=b_1*a+b_2* a$,

$6$. $(b_1*b_2)\cdot(a_1*a_2)=a_1*a_2$,

$7$. $(b_1*b_2)\cdot(a*b)=a*b$,

$8$. $a_1*(b\cdot a_2)=a_1*a_2$,

$9$. $b*(b_1\cdot a)=b*a$,

$10$. $\om(b\cdot a)=\om(b)\cdot\om(a)$,

$11$. $\om(a*b)=\om(a)*b=a*\om(b)$,

$12$. $x*y+z*t=z*t+x*y$,

\noindent for each $\om\in\Om_1'$, $*\in\Om_2{}'$, $b$, $b_1$,
$b_2\in B$, $a,a_1,a_2\in A$ and for  $x,y,z,t\in A\cup B$
whenever each side of $12$ has a sense.
\end{Prop}

Note that in the formulation of Proposition 1.1 in \cite{Da} we mean
that the set of identities of the category of groups with operations
contains only identities from $\bE_G$, but it is not mentioned
there. The same concerns to some other statements of \cite{Da}. In
the case where we have a category $\bC$ with the set of identities
$\bE$, the conditions 1-12 of  the above proposition are necessary
conditions. Of course it is possible according to other identities
included in $\bE$ to write down the corresponding conditions for
derived actions which will be necessary and sufficient for the set
of actions to be a set  of derived actions (i.e. for $B\ltimes
A\in\bC$). Denote all these identities  by $\wt{\bE}_G$ and
$\wt{\bE}$ respectively. If the addition is commutative in $\bC$,
then $\wt{\bE}$ (resp. $\wt{\bE}_G$) consists of the same kind of
identities that we have  in $\bE$ (resp. in $\bE_G$), written down
for the elements from the set $A\cup B$, whenever each identity has
a sense. We will denote by $\wt{\text{Axiom 2}}$  the identities for
the action in $\bC$, which correspond to Axiom 2 (see \cite{Da}). In
the category of groups, Lie, associative, Leibniz algebras derived
actions are called simply actions. We will use this terminology in
these special cases; we will also say ``an action in $\bC$'', if it
is a derived action, and we will say a set of actions is not an
action in $\bC$ if this set is not a set of derived actions. Recall
that a left action of a group $B$ on $A$ is a  map
$\varepsilon:B\times A\lra A$, which we denote by
$\varepsilon(b,a)=b\cdot a$, with conditions
\begin{align*}
(b_1+b_2)\cdot a&={b_1}\cdot ({b_2}\cdot a),\\
0\cdot a&=a,\\
b\cdot (a_1+a_2)&=b\cdot {a_1}+ b\cdot {a_2}.
\end{align*}
The right action  is defined in an analogous way.

All algebras below are considered over a commutative ring {\it k}
with unit.

  In the case of associative algebras an action of $B$ on $A$
is a pair of bilinear maps
\begin{equation}
B\times A\lra A,\quad A\times B\lra A
\end{equation}
which we denote respectively as $(b,a)\mapsto b*a$, $(a,b)\mapsto
a*b$, with conditions
\begin{align*}
&\quad(b_1*b_2)*a=b_1*(b_2*a),\\
&\quad a*(b_1*b_2)=(a*b_1)*b_2,\\
&\quad(b_1*a)*b_2=b_1*(a*b_2),\\
&\quad b*(a_1*a_2)=(b*a_1)*a_2,\\
&\quad(a_1*a_2)*b=a_1*(a_2*b),\\
&\quad a_1*(b*a_2)=(a_1*b)*a_2.\\
\end{align*}
Here the associative algebra operation is denoted by $*$ (resp.
$a_1*a_2$) and the corresponding action by the same sign $*$
(resp. $b*a$).

For Lie algebras an action of $B$ on $A$ is a bilinear map
$B\times A\lra A$, which we denote by $(b,a)\mapsto [b,a]$, with
conditions
\begin{align*}
&\quad[[b_1,b_2],a]=[b_1,[b_2,a]]-[b_2,[b_1,a]],\\
 &\quad[b,[a_1,a_2]]=[a_1,[b,a_2]]+[[b,a_1],a_2].\\
\end{align*}
Note that we actually have above again two bilinear maps (2.4):
$(b,a)\mapsto [b,a]$, $(a,b)\mapsto [a,b]$ with conditions
\begin{align*}
&\quad[b,a]=-[a,b],\\
&\quad[a,[b_1,b_2]]+[b_1,[b_2,a]]+[b_2,[a,b_1]]=0,\\
&\quad[a_1,[b,a_2]]+[b,[a_2,a_1]]+[a_2,[a_1,b]]=0.\\
\end{align*}

Recall from \cite{Lo} that a \textit{Leibniz algebra} L over a
commutative ring   {\it k} with unit  is a {\it k} -module
equipped with a bilinear map $[-,-]: {\rm L} \times {\rm L}
\rightarrow {\rm L}$ which satisfies the following identity,
called the Leibniz identity:
\begin{equation*}
\lbrack x,[y,z]]=[[x,y],z]-[[x,z],y]
\end{equation*}%
for all $x,y,z\in$ L.

Obviously, when $[x,x]=0$ for all $x\in$ L, the Leibniz bracket is
skew-symmetric, therefore the Leibniz identity comes down to the
Jacobi identity, and a Leibniz algebra is then just a Lie algebra.

For Leibniz algebras an action of $B$ on $A$ is a pair of bilinear
maps (2.4), which we denote again by $(b,a)\mapsto[b,a]$,
$(a,b)\mapsto[a,b]$ with conditions: \allowdisplaybreaks
\begin{align*}
[a_1,[a_2,b]]&=[[a_1,a_2],b]-[[a_1,b],a_2],\\
[a_1,[b,a_2]]&=[[a_1,b],a_2]-[[a_1,a_2],b],\\
[b,[a_1,a_2]]&=[[b,a_1],a_2]-[[b,a_2],a_1],\\
[a,[b_1,b_2]]&=[[a,b_1],b_2]-[[a,b_2],b_1],\\
[b_1,[a,b_2]]&=[[b_1,a],b_2]-[[b_1,b_2],a],\\
[b_1,[b_2,a]]&=[[b_1,b_2],a]-[[b_1,a],b_2].
\end{align*}

Recall \cite{Se} that a derivation for a Lie algebra $A$ over a
ring $k$ is a $k$-linear map $D:A\lra A$ with
$$ D[a_1,a_2]=[D(a_1),a_2]+[a_1,D(a_2)]. $$
The set of all derivations Der$(A)$ of $A$, with the operation
defined by
$$  [D,D']=DD'-D'D  $$
is a Lie algebra.

We recall the construction of the $k$-algebra Bim$(A)$ of
bimultipliers of an associative $k$-algebra $A$ \cite{Ho},
\cite{Mc}. An element of Bim$(A)$ is a pair $f=(f*,*f)$ of
$k$-linear maps from $A$ to $A$ with
\begin{align*}
f*(a*a')&=(f*a)*a',\\
(a*a')*f&=a*(a'*f),\\
a*(f*a')&=(a*f)*a'.
\end{align*}
We prefer to use the notation $*f$ instead of $f*^\circ$.  We
denote by $f*a$ (resp. $a*f$) the value $f*(a)$ (resp. $*f(a)$).
Bim$(A)$ is a $k$-module in an obvious way. The operation in
Bim$(A)$ is defined by
$$  f*f'=(f*f'*, *f*f'),    $$
and Bim$(A)$ becomes  a $k$-algebra. Note that here we use
different notations than in \cite{Mc} and \cite{LL}. Here as above
$*$ denotes an operation in associative algebra, and $f*f'*,*f*f'$
denote the compositions of maps. Thus
\begin{align*}
&(f*f'*)(a)=f*(f'*a),\\
&(*f*f')(a)=(a*f)*f'.
\end{align*}
For the addition we have
$$f+f'=((f*)+f'*,*f+(*f')),$$
where
\begin{align*}
&((f*)+f'*)(a)=f*a+f'*a,\\
&(*f+(*f'))(a)=a*f+a*f'.
\end{align*}
For a Leibniz $k$-algebra $A$ we define the $k$-algebra Bider$(A)$
of biderivations in the following way. An element of Bider$(A)$ is
a pair $\vf=([\;,\vf],[\vf,\;])$ of $k$-linear maps $A\lra A$ with
\begin{align*}
\big[[a_1,a_2],\vf\big]&=\big[a_1,[a_2,\vf]\big]+\big[[a_1,\vf],a_2\big],\\
\big[\vf,[a_1,a_2]\big]&=\big[[\vf,a_1],a_2\big]-\big[[\vf,a_2],a_1\big],\\
\big[a_1,[a_2,\vf]\big]&=-\big[a_1,[\vf,a_2]\big].
\end{align*}
We used above the notation:
$[\vf,\;](a)=[\vf,a],[\;,\vf](a)=[a,\vf]$. Biderivations where
defined by Loday in \cite{Lo}, where another notation is used;
biderivation is a pair $(d,D)$, where according to our definition
$[\vf,\;]=D$, $[\;,\vf]=-d$ and instead of the third condition we
have in \cite{Lo} $[a_1,d(a_2)]=[a_1,D(a_2)]$.

The operation in Bider$(A)$ is defined by:
$$ [\vf,\vf']=\big([\;,[\vf,\vf']],[[\vf,\vf'],\;]\big),     $$
where
\begin{align}
\big[a,[\vf,\vf']\big]&=\big[[a,\vf],\vf'\big]-\big[[a,\vf'],\vf\big],\tag{2.5.1}\\
\big[[\vf,\vf'],a\big]&=\big[\vf,[\vf',a]\big]+\big[[\vf,a],\vf'\big].\tag{2.5.2}
\end{align}

Note that we could define $[[\vf,\vf'],\;]$ by
\begin{equation}
\big[[\vf,\vf'],a\big]=-\big[\vf,[a,\vf']\big]+\big[[\vf,a],\vf'\big].\tag{2.5.2'}
\end{equation}
To avoid confusions we forget about $*^\circ$ in special cases, e.g.
for the [\;,\;] operation. The above given both operations define a
Leibniz algebra structure on Bider$(A)$. It is easy to see that the
second definition (2.5.1), (2.5.2') gives the algebra which is
isomorphic to the biderivation algebra defined in \cite{Lo};
according to this definition $[(d,D),(d',D')]=(dd'-d'd,Dd'-d'D)$.

We have a set of actions of Der$(A)$, Bim$(A)$ and Bider$(A)$ on
$A$. These actions are defined by
\begin{align*}
[D,a]&=D(a),\\
f*a&=f*(a),\\
a*f&=*f(a),\\
[\vf,a]&=[\vf,\;](a),\quad [a,\vf]=[\;,\vf](a),
\end{align*}
where $a\in A$, $D\in\text{Der}(A)$, $f=(f*,*f)\in\text{Bim}(A)$,
$\vf=([\;,\vf],[\vf,\;])\in\text{Bider}(A)$ and $A$ is a Lie
algebra, an associative algebra and a Leibniz algebra
respectively.

In the case of Lie algebras the action of Der$(A)$ on $A$ is a set
of derived actions, thus this action satisfies the corresponding
conditions of an action in $\bL$ie, but for the case of
associative and Leibniz algebras these actions do not satisfy all
the conditions given above respectively for the action in
$\mathbb{A}$ss and $\bL$eib. Note that for the case of Leibniz
algebras if $[\vf,[\vf',a]]=-[\vf,[a,\vf']]$ for any $a\in A$ and
$\vf,\vf'\in\text{Bider}(A)$, then above two ways of defining
operations in Bider$(A)$ are equal and also the action of
Bider$(A)$ becomes a derived action (see below Proposition ~4.8).

We have an analogous situation for associative algebras. The
action of Bim$(A)$ on $A$ is not a derived action because the
condition \addtocounter{equation}{+1}
\begin{equation}
(f*a)*f'=f*(a*f')
\end{equation}
fails. So if we would have the condition for associative algebra
$A$ that for any two bimultipliers is fulfilled (2.6), then the
action of Bim$(A)$ on $A$ defined above is a set of derived
actions on $A$ (see below Proposition 4.7).

An {\it alternative algebra} $A$ over a field $F$ is an algebra
which satisfies the identities $$x^2y = x(xy)$$ and $$yx^2 =
(yx)x$$ for all $x,y \in A$. These identities are known
respectively as the left and right alternative laws. We denote the
corresponding category of alternative algebras by $\mathbb{A}$lt.
Clearly any associative algebra is alternative. The class of
8-dimensional Cayley algebras is an important class of alternative
algebras which are not associative \cite{Sc}.

The axioms above for alternative algebras are equivalent to the
following:
$$ x(yz) = (xy)z+(yx)z-y(xz)$$
and
$$(xy)z =x(yz) - (xz)y + x(zy)$$
We consider these conditions as Axiom 2 and consequently alternative
algebras can be interpreted as categories of interest.

For alternative algebras over a field $F$ an action of $B$ on $A$
is a pair of bilinear maps (2.4), which we denote again by
$(b,a)\mapsto b a$, $(a,b)\mapsto a b$ with conditions:
\allowdisplaybreaks
\begin{align*}
b (a_1 a_2)& = (b a_1) a_2 + (a_1 b)a_2-a_1(b a_2), \\
(a_1 a_2) b & = a_1 (a_2 b) - (a_1 b) a_2 + a_1 (b a_2), \\
(b a_1) a_2 & = b (a_1 a_2) - (b a_2) a_1 + b (a_2 a_1), \\
a_1(a_2 b)&=(a_1 a_2)b + (a_2 a_1) b - a_2(a_1 b),\\
(b_1 b_2) a & = b_1 (b_2 a) - (b_1 a) b_2 + b_1 (a b_2), \\
a (b_1 b_2) &= (a b_1) b_2 + (b_1 a) b_2 - b_1(a b_2), \\
(a b_1) b_2& = a(b_1 b_2) - (a b_2) b_1 + a(b_2 b_1), \\
b_1(b_2a)&=(b_1 b_2) a+ (b_2 b_1) a - b_2(b_1 a). \\
\end{align*}

A crossed module in $\bC$ is a triple $(C_0,C_1,\pa)$, where
$C_0,C_1\in\bC$, $C_0$ acts on $C_1$ (i.e. we have a derived
action in $\bC$) and $\pa:C_1\lra C_0$ is a morphism in $\bC$ with
conditions:

(i) $\pa(r\cdot c)=r+\pa(c)-r$;

(ii) $\pa(c)\cdot c'=c+c'-c$;

(iii) $\pa(c)*c'=c*c'$;

(iv) $\pa(r*c)=r*\pa(c)$, $\pa(c*r)=\pa(c)*r$

\noindent for any $r\in C_0$, $c,c'\in C_1$, and $*\in \Om_2{}'$.

A morphism between two crossed modules
$(C_0,C_1,\pa)\lra(C'_0,C'_1,\pa')$ is a pair of morphisms
$(T_0,T_1)$ in $\bC$, $T_0:C_0\lra C_0'$, $T_1:C_1\lra C_1'$, such
that
\begin{align*}
T_0\pa(c)&=\pa'T_1(c),\\
T_1(r\cdot c)&=T_0(r)\cdot T_1(c),\\
T_1(r*c)&=T_0(r)*T_1(c)
\end{align*}
for any $r\in C_0$, $c\in C_1$ and $*\in\Om_2{}'$.

\begin{defi}
For any object $A$ in $\bC$ an actor of $A$ is a crossed module
$\pa:A\lra \text{\rm Actor}(A)$, such that for any object $C$ of
$\bC$ and an action of $C$ on $A$ there is a unique morphism
$\vf:C\lra\text{\rm Actor}(A)$ with $c\cdot a=\vf(c)\cdot a$,
$c*a=\vf(c)*a$ for any $*\in\Om_2{}'$, $a\in A$ and $c\in C$.
\end{defi}

See the equivalent Definition 3.9. in Section 3.

From this definition it follows that an actor object
$\text{Actor}(A)$, for the object $A\in\bC$, with this properties is
a unique object up to an isomorphism in $\bC$.

Note that according to the universal property of an actor object,
for any two elements $x,y$ in Actor$(A)$ from
$x\os{\cdot}{*}a=y\os{\cdot}{*}a$, (here we mean equalities for
the dot action and the action $*$, for any $*\in\Om_2{}'$ and any
$a \in A$) and $(w_1 \cdots w_n x)\cdot a=(w_1 \cdots w_n y)\cdot
a, w_1 \cdots w_n \in \Omega'_1,$ it follows that  $x=y$.

It is well-known that for the case of groups
Actor$(G)=\text{Aut}(G)$; the corresponding crossed module is
$\pa:G\lra\text{Aut}(G)$, where $\pa$ sends any $g\in G$ to the
inner automorphism of $G$ defined by $g$ (i.e.
$\pa(g)(g')=g+g'-g$, $g'\in G$). For the case of Lie algebras
Actor$(A)=\text{Der}(A)$, $A\in \bL$ie, and the operator
homomorphism $\pa:A\lra\text{Der}(A)$ is defined by
$\pa(a)=[a,\;]$, so $\pa(a)(a')=[a,a']$.

As we have seen above, in general, in $\mathbb{A}$ss and $\bL$eib
the objects Bim$(A)$ and Bider$(A)$ do not have derived actions on
$A$ in the corresponding categories. So the obvious homomorphisms
$A\lra\text{Bim}(A)$, $A\lra\text{Bider}(A)$ do not define crossed
modules in $\mathbb{A}$ss and $\bL$eib for any $A$ from $
\mathbb{A}$ss and $\bL$eib respectively.

It is well-known \cite{No} that for the case of groups if $N$ is a
normal subgroup of $G$ and $\tau:N\lra \text{Inn}(N)$ is the
homomorphism sending any element $n$ to the corresponding inner
automorphism $(\tau(n)(n')=n+n'-n)$, since $G$ acts on $N$ by
conjugation, we have a unique homomorphism $\tht:G\lra
\text{Actor}(N)$, with $\tht(g)\cdot n=g\cdot n$. Inn$(N)$ is a
normal subgroup of Actor$(N)$, $\tht$ extends $\tau$ and we have a
commutative diagram
\begin{gather}
\xymatrix{0\ar[r]&N\ar[r]\ar[d]_-{\tau}&G\ar[r]\ar[d]_-{\tht}&G/N\ar[r]\ar@{-->}[d]&0\\
0\ar[r]&\text{Inn}(N)\ar[r]&\text{Actor}(N)\ar[r]&\text{Out}(N)\ar[r]&0.}   \notag\\[-12mm]
 \
\end{gather}
\

According to the work of R. Lavendhomme and Th. Lucas \cite{LL} in
the categories $\bG r$, $\bL$ie the actor crossed modules $A\lra
\text{Actor}(A)$ are terminal objects in the categories of crossed
modules under $A$. If Ann$(A)=(0)$ or $A^2=A$ then Bim$(A)$ acts on
$A$ and the corresponding crossed module $A\lra\text{Bim}(A)$ is a
terminal object in the category of crossed modules under $A$. It is
easy to see that in this case
$$    \text{Bim}(A)=\text{Actor}(A).  $$

\begin{defi}
A general actor object {\rm GActor}$(A)$ for $A$, $A\in\bC$, is an
object from $\bC_G$, which has a set of actions on $A$, which is a
set of derived actions in $\bC_G$, i.e. satisfies conditions of
Proposition 2.6, there is a morphism $d:A\lra $ {\rm GActor}$(A)$
in $\bC_G$ which defines a crossed module in $\bC_G$ and for any
object $C\in\bC$ and a derived action of $C$ on $A$, there exists
in $\bC_G$ a unique morphism $\vf:C\lra$ {\rm GActor}$(A)$ such
that $c\os{\cdot}{*}a=\vf(c)\os{\cdot}{*}a$,  for any $c\in C$,
$a\in A$, $*\in\Om_2{}'$.
\end{defi}

It is easy to see that Bim$(A)$ and Bider$(A)$ are general actor
objects for $A\in \mathbb{A}$ss, $A\in\bL$eib respectively. These
constructions satisfy the existence of the commutative diagram
like (2.7).

\section{The  main construction}

In this section $\bC$ will denote a category of interest with a set
of operations $\Om$ and with a set of identities $\bE$. Let $\bC_G$
be the corresponding general category of groups with operations.
According  to the definition given in Section 2, $\bC_G$ is a
$\{\Om,\bE_G\}$-algebra, where $\Om$ is the same set of operations
as we have in $\bC$, and $\bE_G$ includes group laws and identities
(c) and (d) from the definition of Section 2. Thus we have
$\bE_G\hookrightarrow\bE$ and $\bC$ is the full subcategory  of
$\bC_G$, $\bC\hookrightarrow\bC_G$. Let $A\in\bC$; consider all
split extensions of $A$ in $\bC$
$$  \xymatrix{E_j:0\ar[r]&A\ar[r]^-{i_j}&C_j\ar[r]^-{p_j}&B_j\ar[r]&0},\quad j\in\bJ.   $$
Note that it may happen that $B_j=B_k=B$, for $j\neq k$, then these
extensions will correspondent to different actions of $B$ on $A$.
 Let $\{b_j\cdot,b_j*|b_j\in
B_j,\;*\in\Om_2'\}$ be the corresponding set of derived actions for
$j\in\bJ$. For any element $b_j\in B_j$ denote
${\bbb}_j=\{b_j\cdot,b_j*,\;*\in\Om_2'\}$. Let
$\bB=\{{\bbb}_j|b_j\in B_j,\;j\in\bJ\}$.

Thus each element ${\bbb}_j\in\bB$, $j\in\bJ$ is a special type of
function ${\bbb}_j:\Om_2\lra\text{Maps}(A\to A)$,
${\bbb}_j(*)=b_j*-:A\lra A$.

According to Axiom 2 from the definition of a category of interest,
we define $*$ operation, ${\bbb}_i*{\bbb}_k,*\in\Om_2'$, for the
elements of $\bB$ by the equalities
\begin{align*}
&({\bbb}_i*{\bbb}_k)\ol{*}\,(a)=W(b_i,b_k;a;*,\ol{*}),\\
&({\bbb}_i*{\bbb}_k)\cdot(a)=a.
\end{align*}
We define the operation of addition by
\begin{align*}
&({\bbb}_i+{\bbb}_k)\cdot(a)=b_i\cdot(b_k\cdot a),\\
&({\bbb}_i+{\bbb}_k)*(a)=b_i*a+b_k*a.
\end{align*}
For a unary operation $\om\in\Om_1'$ we define
\begin{align*}
&\om({\bbb}_k)\cdot (a)=\om(b_k)\cdot (a),\\
&\om({\bbb}_k)* (a)=\om(b_k)* (a),\\
&\om(b*b')=\om(b)*b'\textrm{ and we will have }\om(b)*b'=b*\om(b'),\\
&\om(b_1+\cdots+b_n)=\om(b_1)+\cdots+\om(b_n),\\
&(-{\bbb}_k)\cdot (a)=(-b_k)\cdot a,\\
&(-b)\cdot (a)=a\\
&(-{\bbb}_k)*(a)=-(b_k * a),\\
&(-b)*(a)=-(b*(a)),\\
&- (b_1+\cdots+b_n)=-b_n-\cdots-b_1,\\
\end{align*}
where $b,b',b_1,...,b_n$ are certain combinations of star
operations on the elements of $\bB$, i.e. the elements of the type
$\bbb_{i_1}*_1\cdots*_{n-1}\bbb_{i_n}$, $n>1$.

We do not know if the new functions defined by us are again in
$\bB$. Denote by $\cB(A)$ the set of functions
($\Om_2\lra\text{Maps}(A\to A))$ obtained by performing all kind of
above defined operations on elements of $\bB$ and new obtained
elements as the results of operations. Note that $b=b'$ in $\cB(A)$
means that $b\os{\cdot}{*}a=b'\os{\cdot}{*}a$, $w_1 \dots w_n b
\cdot a = w_1 \dots w_n b' \cdot a$ for any $a\in A$, $*\in\Om_2'$,
$w_1 \dots w_n \in \Omega'_1$ and for any $n$. It is an equivalence
relation and under $\cB(A)$ we mean the corresponding quotient
object.

\begin{Prop}
$\cB(A)$ is an object of  $\bC_G$.
\end{Prop}

\begin{proof}
Direct easy checking of the identities.
\end{proof}

As above, we will write for simplicity $b\cdot(a)$ and $b*(a)$
instead of $(b(+))(a)$ and $(b(*))(a)$ for $b\in \cB(A)$ and $a\in
A$. Define the set of actions of $\cB(A)$ on $A$ in a natural way.
For $b\in \cB(A)$ we define $b\cdot a=b\cdot(a)$, $b*a=b*(a)$,
$*\in\Om_2'$. Thus if $b=\bbb_{i_1}*_1\cdots*_{n-1}\bbb_{i_n}$,
where we mean certain brackets, we have
\begin{align*}
&b\,\ol{*}\,a=(\bbb_{i_1}*_1\cdots*_{n-1}\bbb_{i_n})\,\ol{*}\,(a),\\
&b\,\cdot\, a=a.
\end{align*}
The right side of the equality is defined inductively according to
Axiom 2. For $b_k\in B_k$, $k\in\bJ$, we have
\begin{align*}
\bbb_k*a=\bbb_k*(a)=b_k*a,\\
\bbb_k\cdot a=\bbb_k\cdot(a)=b_k\cdot a.
\end{align*}
Also
\begin{align*}
&(b_1+b_2+\dots+b_n)*a=b_1*(a)+\dots+b_n*(a),\\
&\hskip+6.5cm\text{for}\quad b_i\in \cB(A),\quad i=1, \cdots, n\\
&(b_1+b_2+\dots+b_n)\cdot a=b_1\cdot(b_2\cdots(b_n\cdot(a))\cdots),\\
&\hskip+6.5cm b_i\in \cB(A),\quad i=1, \cdots,n\\
&\om(b)\cdot a=a\quad\text{if}\quad b=b_1*\cdots*b_n,\quad b_i\in \cB(A),\quad i=1,\cdots,n\\
&\om(\bbb_k)\cdot a=\om(b_k)\cdot a,\quad k\in\bJ,\quad b_k\in
B_k.
\end{align*}

\begin{Prop}
The set of actions of $\cB(A)$ on $A$ is a set of derived actions
in $\bC_G$.
\end{Prop}

\begin{proof} The checking shows that  the set of  actions of $\cB(A)$ on $A$ satisfies conditions of Proposition 2.6,
 which proves that it is a set of derived actions in $\bC_G$.
\end{proof}

Define the map $d:A\lra \cB(A)$ by $d(a)=\aaa$, where
$\aaa=\{a\cdot,a*,*\in \Om_2'\}$. Thus we have by definition
\begin{align*}
&d(a)\cdot a'=a+a'-a,\\
&d(a)*a'=a*a',\quad \forall a,a'\in A,\quad *\in\Om_2'.
\end{align*}

\begin{Lem}
$d$ is a homomorphism in $\bC_G$.
\end{Lem}

\begin{proof}
We have to show that $d(\om a)=\om d(a)$ for any $\om\in\Om_1'$. For
this we need to show that
\begin{align*}
&d(\om a)\cdot(a')=(\om d(a))\cdot(a')\quad \\
&\om'(d(\om a))\cdot a' = \om'(\om d(a)) \cdot a', \quad \text{for any $\om' \in \Omega'_1$}\\
 &d(\om
a)*(a')=(\om d(a))*(a'),\quad\text{for any}\quad *\in\Om_2'.
\end{align*}
We have
\begin{align*}
d(\om a)\cdot a'&=\om a+a'-\om a,\\
\om d(a)\cdot a'&=\om(\aaa)\cdot a'=\om a+a'-\om a,
\end{align*}
The second equality follows form the first one. For the third
equality we have
\begin{align*}
d(\om a)* a'&=(\om a)* a',\\
(\om d(a))* a'&=\om(\aaa)* a'=\om(a)*a'
\end{align*}
for $\om=-$ we have to show $d(-a)\cdot(a')=(-(da))\cdot a'$ and
$(d(-a))*a'=(-d(a))*a'$. The checking of these equalities is an
easy exercise.

Now we will show  that $d(a_1+a_2)=d(a_1)+d(a_2)$. The direct
computation   of both sides for each $a\in A$ gives
\begin{align*}
&d(a_1+a_2)\cdot (a)=a_1+a_2+a-a_2-a_1,\\
&\big(d(a_1)+d(a_2)\big)\cdot(a)=d(a_1)\cdot\big(d(a_2)\cdot
a\big),
\end{align*}
which shows that the desired equality holds for the dot action.
The proof of $\om(d(a_1+a_2)) \cdot a = \om(d(a_1)+d(a_2))\cdot a$
is based on the first equality, the property of unary operations
to respect the addition and the fact that $d$ commutes with unary
operations.

For any $*\in\Om_2'$ we shall show that
$$d(a_1+a_2)*(a)=\big(d(a_1)+d(a_2)\big)*(a).$$
We have
\begin{align*}
&d(a_1+a_2)*(a)=(a_1+a_2)*a=a_1*a+a_2*a,\\
&\big(d(a_1)+d(a_2)\big)*(a)=d(a_1)*a+d(a_2)*a=a_1*a+a_2*a
\end{align*}
which proves the equality.

The next equality we have to prove is $d(a_1*a_2)=d(a_1)*d(a_2)$.
For this we need to show that
$d(a_1*a_2)\cdot(a)=(d(a_1)*d(a_2))\cdot(a)$, $\om(d(a_1 * a_2))
\cdot a = \om(d(a_1)*d(a_2))\cdot a$ and
$d(a_1*a_2)\ol{*}(a)=(d(a_1)*d(a_2))\ol{*}(a)$, for any
$\ol{*}\in\Om_2'$.

We have $d(a_1*a_2)\cdot a=a_1*a_2+a-a_1*a_2=a$, since $A\in\bC$
and  therefore it satisfies Axiom 1.

$(d(a_1)*d(a_2))\cdot a=a$, by the definition of the action of
$\cB(A)$ on $A$. The next equality is proved in a similar way
applying that $d$ commutes with $\om$ and $\om(a_1*a_2) = \om(a_1)
*a_2$.

For the next above given identity we have the following
computations:
\begin{align*}
& d(a_1*a_2)\ol{*}(a)=(a_1*a_2)\ol{*}a=W(a_1,a_2;a;*,\ol{*}),\\
&\big(d(a_1)*d(a_2)\big)\ol{*}(a)=W\big(d(a_1),d(a_2);a;*,\ol{*}\big).
\end{align*}

These two expressions on the right sides of above equalities are
equal, by the type of the word $W$ in Axiom 2 and the definition of
$d$.
\end{proof}

\begin{Prop}
$d:A\lra \cB(A)$ is crossed module in $\bC_G$.
\end{Prop}

\begin{proof}
We have to check conditions (i)-(iv) from the definition of a
crossed module  given in Section 2.

(i) condition states that $d(b\cdot a)=b+d(a)-b$ for $a\in A$,
$b\in \cB(A)$; so we have to show that $d(b\cdot
a)\os{\cdot}{*}a'=(b+da-b)\os{\cdot}{*}a'$ and $\om_1 \dots
\om_n(d(b\cdot a)) \cdot a' = \om_1 \dots \om_n(b+da-b)\cdot a'$.
Below we compute each side for the dot action of the first
equality:
\begin{align*}
&d(b\cdot a)\cdot a'=b\cdot a+a'-b\cdot a,\\
&(b+d(a)-b)\cdot a'=b\cdot(d(a)\cdot(-b\cdot a'))\\
&\quad=b\cdot(a-b\cdot a'-a)=b\cdot a+a'-b\cdot a.
\end{align*}
The second equality is proved in a similar way. Now we compute
each side of the first equality for the $*$ action. $d(b\cdot
a)*a'=(b\cdot a)*a'=a*a'$ by Proposition 2.6;
$(b+da-b)*a'=b*a'+d(a)*a'-b*a'=b*a'+a*a'-b*a'=a*a'$, here we
applied Axiom 1, that $\ol{a}+a*a'=a*a'+\ol{a}$, for any element
$\ol{a}$ of $A$.

 We have to show: (ii) $d(a_1)\cdot
a_2=a_1+a_2-a_1$, (iii) $d(a_1)*a_2=a_1*a_2$; both (ii) and (iii)
are true by definition of $d$. Note that
$a_1*(d(a_2))=d(a_2)*^\circ a_1=a_2*^\circ a_1=a_1*a_2$.

The first condition of (iv) states
$$    d(b*a)=b*d(a),\quad\text{for any}\quad b\in \cB(A),\quad a\in A, *\in\Om_2'.$$
Thus we have to show
\begin{equation}
d(b*a)\os{\cdot}{\ol{*}}\,a'=(b*d(a))\os{\cdot}{\ol{*}}\,a', \
\om(db*a)\cdot a' = \om(b*da)\cdot a', \quad\text{for}\quad  \om
\in \Omega'_1, \ol{*}\in\Om_2'.
\end{equation}
First we show (3.1) for the dot operation. The second equality for
the dot operation is proved similarly applying properties of unary
operations. The right side of (3.1) in this case is equal to $a'$.
For the left side we obtain
$$    d(b*a)\cdot a'=b*a+a'-b*a.  $$
If $b=\bbb_i$, then $b*a=b_i*a$ and since $B_i\in\bC$, and $B_i$
acts on $A$ (action is in $\bC$), by Axiom 1 for the action of
$B_j$ on $A$ we shall have $b*a+a'=a'+b*a$ and so $d(b*a)\cdot
a'=a'$.

If $b=\bbb_{i_1}*_1\cdots*_{n-1}\bbb_{i_n}$ then, by the
definition of $*$ operation in $\cB(A)$, $b*a$ is the sum of the
elements of the type $b_{i_t}*\ol{a}_t$ for certain $i_t$ and the
element $\ol{a}_t\in A$; this kind of element again commutes with
any element of $A$. So that $d(b*a)\cdot a'=a'$. We will have the
same result if $b$ is the sum of the elements  of the type
$\bbb_{i_1}*_1\cdots*_{n-1}\bbb_{i_n}$.

Now we shall show (3.1) for the $*$ operation. By the definition
of $d$ we have
$$    d(b*a)\ol{*}\,a'=(b*a)\ol{*}\,a'.  $$
In the case $b=\bbb_i$, $i\in\bJ$, $b*a=\bbb_i*a=b_i*a$, so we
obtain
$$    d(b*a)\ol{*}\,a'=(b_i*a)\ol{*}\,a=W(b_i,a;a';*,\ol{*}).  $$
We have the last equality according to the properties of an action
in $\bC$, which correspond to Axiom 2. For the right side of (3.1)
in case $b=\bbb_i$ we have
$$(b*d(a))\ol{*}\,a'=(\bbb_i*\aaa)\ol*\,a'=W(b_i,a;a';*,\ol{*}).$$
Suppose $b=\bbb_{i_1}*_1\cdots *_{n-1}\bbb_{i_n}$, then in the same
way as it was in the previous proof, $b*a$ is the sum of the
elements of the type $b_{i_t}*\ol{a}_t$ and $(b*a)\ol{*}\,a'$ is the
sum of the elements of the type $(b_{i_t}*\ol{a}_t)\ol{*}\,a'$. The
element from the right side of (3.1) will be the same type of the
sum of the elements $(b_{i_t}*\ol{a}_t)*a'$. Applying
$\wt{\text{Axiom 2}}$ to the element $(b_{i_t}*\ol{a}_t)\ol{*}\,a'$,
by the definition of the operation for the elements of $\cB(A)$ (for
the element $(\bbb_{i_t}*\ol{\aaa}_t)*a')$ and from the facts that
$\bbb_{i_t}*a=b_{i_t}*a$, $\ol{\aaa}_t*a=\ol{a}_t*a$, we will have
the desired equality (3.1). In the analogous way we will prove (3.1)
for $*$ operation in case $b$ is a sum of the elements of the form
$\bbb_{i_1}*_1\cdots *_{n-1}\bbb_{i_n}$. The second condition of
(iv) can be proved in a similar way.
\end{proof}

\begin{Prop}
If $A$ has an actor in $\bC$, then $\cB(A)={\rm Actor}(A)$.
\end{Prop}

\begin{proof}
From the existence of Actor$(A)$ it follows that Actor$(A)$ is one
of the objects $B_i$, which acts on $A$. We have a natural
homomorphism $e:\text{Actor}(A)\lra \cB(A)$ in $\bC_G$ sending $b_i$
to $\bbb_i$, $b_i\in B_i$.  According to the note made in Section 2,
if $b_i\neq b_i'$ in Actor$(A)$, then $\bbb_i\neq\bbb_i'$; thus $e$
is an injective homomorphism. Let $\vf_j:B_j\lra\text{Actor}(A)$ be
a unique morphism with $\vf(b_j)\os{\cdot}{*}a=b_j\os{\cdot}{*}a$,
$b_j\in B_j$, $j\in\bJ$, $a\in A$. $e$ is a surjective homomorphism,
since for any element $\bbb_{i_1}*_1\cdots *_{n-1}\bbb_{i_n}$ of
$\cB(A)$ there exists the element $\vf_{i_1}(b_{i_1})*_1\cdots
*_{n-1}\vf_{i_n}(b_{i_n})$ in Actor$(A)$ with
$e(\vf_{i_1}(b_{i_1})*_1\cdots*_{n-1}\vf_{i_n}(b_{i_n}))=\bbb_{i_1}*_1\cdots
*_{n-1}\bbb_n$ which ends the proof.
\end{proof}

\begin{theo}
Let $\bC$ be a category of interest and $A\in\bC$. $A$ has an
actor if and only if $\cB(A)\ltimes A\in\bC$. If it is the case,
then {\rm Actor}$(A)=\cB(A)$.
\end{theo}

\begin{proof}
From the Proposition 3.5 it follows that if $A$ has an actor then
$\cB(A)\in\bC$ and $\cB(A)$ has a derived action on $A$. By the
theorem of Orzech \cite{Or} (see Section 2,  Theorem 2.5) we will
have $\cB(A)\ltimes A\in\bC$. The converse is also easy to prove.
Since $\cB(A)\ltimes A\in\bC$, from the split exact sequence
$\xymatrix{0\ar[r]&A\ar[r]^-{i}&\cB(A)\ltimes A \ar[r]&
\cB(A)\ar[r]&0}$ in $\bC_G$, $\cB(A) $ = Coker $i$ and thus it is an
object of $\bC$; again by Theorem 2.5 $\cB(A)$ has a derived action
on $A$ in $\bC$ (it is the action we have defined). By Proposition
3.4, $d:A\lra \cB(A)$ is a crossed module in $\bC_G$; since
$\cB(A)\in\bC$, and the action of $\cB(A)$ on $A$ is a derived
action in $\bC$, it follows that $d:A\lra \cB(A)$ is a crossed
module in $\bC$. Now we have to show the universal property of this
crossed module. For any action of $B_k$ on $A$, $k\in\bJ$, we define
$\vf_k:B_k\lra \cB(A)$ by $\vf_k(b_k)=\bbb_k$, for any $b_k\in B_k$,
where $\bbb_k\in\bB$. By definition of $\bB$,
$\bbb_k\os{\cdot}{*}a=b_k\os{\cdot}{*}a$, $*\in\Om_2'$, and we
obtain $\vf_k(b_k)\os{\cdot}{*} a=b_k\os{\cdot}{*} a$; $\vf_k$ is a
homomorphism in $\bC$. For another homomorphism $\vf'_k$ we would
have
$\vf'_k(b_k)\os{\cdot}{*}a=b_k\os{\cdot}{*}a=\vf_k(b_k)\os{\cdot}{*}a$,
$\om(\vf'_k(b_k)) \cdot a = \vf'_k(\om b_k)\cdot a = (\om b_k)\cdot
a = \om(\vf(b_k))\cdot a$, for any $b_k\in B_k$, $a\in A$, $\om \in
\Omega'_1$, and $*\in\Om_2'$, which means that
$\vf_k(b_k)=\vf'_k(b_k)$, for any $b_k\in B_k$, this gives the
equality $\vf_k=\vf'_k$, which proves the theorem.
\end{proof}

\begin{theo}
Let $\bC$ be a category of interest. For any $A\in \bC$, $\cB(A)=$
{\rm GActor({\it A})}.
\end{theo}
\begin{proof}
By Propositions 3.2 and 3.4 and Lemma 3.3 we have the crossed module
$d : A \to \cB(A)$  in $\mathbb{C}_G$. For any object $C \in
\mathbb{C}$ which has a  derived action on $A$ we construct the
homomorphism $\varphi : C \to \cB(A)$ in $\mathbb{C}_G$ with the
property $c \os{\cdot}{*} a = \varphi(c) \os{\cdot}{*} a$ and show
that $\varphi$ is unique with this property in the similar way as we
have done for $\varphi_k$ in the proof of Theorem 3.6.
\end{proof}

Below we give a categorical presentation of the necessary and
sufficient conditions for the existence of an actor in the
category of interest $\mathbb{C}$, i.e. for any object $A \in
\mathbb{C}$. Actually we have constructed the  functor $T : \cB(-)
\ltimes (-) : \mathbb{C} \to \mathbb{C}_G$. This functor is
defined in a natural way $T(A) = \cB(A) \ltimes A, A \in
\mathbb{C}$. For the definition of $T(\alpha), \alpha:A \to A'$ in
$\mathbb{C}$, we apply the universality property of the general
actor object in the following way. The pushout diagram in
$\mathbb{C}_G$
$$   \xymatrix{0\ar[r]&A\ar[r]^-{i}\ar[d]_-{\alpha}&{\cB(A)}\ltimes A\ar[r]^-{p}\ar[d]
&{\cB(A)} \ar[r] \ar@{=}[d]&0\\
0\ar[r]&A'\ar[r]&C\ar[r]&{\cB(A)}\ar[r]&0,}
$$
where the first sequence is split, implies that the second one is
also split. Thus  by Theorem 2.5,  ${\cB(A)}$ has a derived set of
actions on $A$ in $\mathbb{C}_G$. The object ${\cB(A')}$ is a
general actor object for $A'$ in $\mathbb{C}_G$. Thus there exists
a unique arrow $\varphi : {\cB(A)} \to {\cB(A')}$ in
$\mathbb{C}_G$ such that  $\varphi(b)*a=b*a$ and $\varphi(b)\cdot
a=b\cdot a$. We define $T(\alpha) :{\cB(A)}\ltimes A \to
{\cB(A')}\ltimes A'$ by $T(\alpha)(b,a) = (\varphi(b),\alpha(a))$.
It is easy to check that   $T(\alpha)$ is a homomorphism in
$\mathbb{C}_G$.

We denote by $Q : \mathbb{C}_G \to \mathbb{C}$ the functor which
assigns to each object $ C \in \mathbb{C}_G$ the greatest quotient
object of $C$ which belongs to $\mathbb{C}$. The above description
gives to Theorem 3.6 the following form.
\bigskip

\noindent \textbf{{Theorem 3.6 $'$}}  {\it Let $\mathbb{C}$ be a
category of interest. There exists an {\rm Actor$(A)$} for any $A\in
\mathbb{C}$ if and only if the following diagram commutes
$$   \xymatrix{\mathbb{C}\ar[r]^-{T} \ar[d]_-{T} &\mathbb{C}_G \ar[d]_-{Q}\\
\mathbb{C}_G &\mathbb{C}\ar[l]^-{E}}
$$
where $E$ denotes the natural inclusion functor.}

\bigskip

Suppose $I$ is an ideal of $C$ in $\bC$ and Actor$(I)$ exists.
Thus we have the  crossed module $d:I\lra\text{Actor}(I)$. Denote
Im $d$ = Inn$(I)$. Thus we have
$$   \text{Inn}(I)=\{\aaa\in\text{Actor}(I)|a\in I\}.   $$
Recall that by definition of $d$, $d(a)=\aaa$, and $\aaa$ is
defined by
\begin{align*}
\aaa\cdot(a')&=a+a'-a,\\
\aaa*(a')&=a*a'.
\end{align*}
It is easy to see that Inn$(I)$ is an ideal of Actor$(I)$. It
follows from the fact that $d:I\lra\text{Actor}(I)$ is a crossed
module and it can also be checked directly. Since $I$ is an ideal
of $C$, we have an action of $C$ on $I$, defined by $c\cdot
a=c+a-c$, $c*a=c*a$, $*\in\Om_2'$. It is a derived action. Thus
there exists a unique homomorphism $\tht:C\lra\text{Actor}(I)$,
such that
$$   \tht(c)\os{\cdot}{*}a=c \os{\cdot}{*}a,\quad a\in I,\quad c\in C,\quad *\in\Om_2'.   $$
Let $\tau:I\lra\text{Inn}(I)$ be a homomorphism defined by $d$,
then $\tht$ induces the commutative diagram
$$   \xymatrix{0\ar[r]&I\ar[r]\ar[d]_-{\tau}&C\ar[r]\ar[d]_-{\tht}
&C/I\ar[r]\ar@{-->}[d]&0\\
0\ar[r]&\text{Inn}(I)\ar[r]&\text{Actor}(I)\ar[r]&\text{Out}(I)\ar[r]&0,}
$$ which is well-known for the case of groups \cite{No} (see Section 2).

For any object $C \in \mathbb{C}$ there is an action of $A$ on
itself defined by $a \cdot a' = a +a' -a; a * a' = a * a'$, for
any $a,a' \in A, * \in \Omega'_2$, where $*$ on the left side
denotes the action and on the right side the operation in $A$. We
call this action the conjugation.

Let $E_A : \xymatrix{0  \ar[r] & A \ar[r] & A \ltimes A \ar[r] &A
\ar[r] \ar@<-1ex>[l] & 0}$ be the split extension which corresponds
to the action of $A$ on itself by conjugation. Consider the category
of all split extensions with fixed $A$; thus the objects are
\[
\xymatrix{0 \ar[r] & A \ar[r] & C  \ar[r] &C' \ar[r] \ar@<-1ex>[l]
& 0,}
\]
and the arrows are triples $(1_A,\gamma,\gamma')$ between
extensions which commute with section homomorphism too.

\begin{Prop}
If $E_t : \xymatrix{0  \ar[r] & A \ar[r] & C \ar[r] &B \ar[r]
\ar@<-1ex>[l] & 0}$ is a terminal object in the category of split
extensions with fixed $A$, then the unique arrow $(1,\gamma,\beta)
: E_A \to E_t$ defines a crossed module $\beta : A \to B$, which
is an actor of $A$.
\end{Prop}

\begin{proof}
The prove  is similar to that of Proposition 3.4. It is obvious that
$B$ has the universal property of an actor. We have to prove that
$\beta : A \to B$ is a crossed module, thus we shall show the
following identities
\begin{align*}
&\beta(a) \cdot a' = a + a'-a,\\
&\beta(b\cdot a) = b + \beta(a) -b,\\
&\beta(a) * a' = a * a'\\
& \beta(b *a)=b * \beta(a).
\end{align*}
for any $a\in A, b \in B, * \in \Omega'_2$. We have the
commutative diagram $$ \xymatrix{E_A : \quad 0  \ar[r] & A \ar[r]
\ar@{=}[d] & A \ltimes A \ar[r] \ar[d]_{\gamma} &A \ar[r]
\ar[d]_{\beta} \ar@<-1ex>[l] & 0\\ E_t : \quad 0  \ar[r] & A
\ar[r] & C \ar[r] &B \ar[r] \ar@<-1ex>[l] & 0}$$ from which we
obtain $\beta(a) \cdot a' = a + a' -a$ and $\beta(a) * a' = a *
a'$ for any $a, a' \in A, * \in \Omega'_2$, which proves the first
and third equalities. Since $E_t$ is a terminal  extension, it has
the following property: if for $b,b' \in B$ we have $b
\os{\cdot}{*} a = b' \os{\cdot}{*} a, \ \omega_1 \cdots
\omega_n(b) \os{\cdot}{*} a = \omega_1 \cdots \omega_n \os{\cdot}
(b') \os{\cdot}{*} a$ for any $a \in A$ and any unary operations
$\omega_1, \cdots, \omega_n \in \Omega'_1, n \in \mathbb{N}$, then
$b = b'$.

For the second equality we have $$(\beta(b \cdot a)) \cdot a' = b
\cdot a + a' - b \cdot a$$
$$(b + \beta(a) -b) \cdot a' = b \cdot(\beta(a) \cdot(-b \cdot
a')) = b \cdot(a-b \cdot a' -a) = b\cdot a + a' - b \cdot a$$
$$(\beta(b \cdot a) )* a'=(b \cdot a) * a' = a * a'$$
by condition 8 of Proposition 2.6.

For the forth equality we have: $$\beta(b *a) \cdot a' = (b *a)
\cdot a' = a'$$ it follows from the property of the derived action
in categories of interest, as a result of Axiom 1 \cite{Da}. The
same property gives

$$(b * \beta(a)) \cdot a' = a'.$$

For a star operation we have:

$$\beta(b *a)*a' = (b *a)*a'$$
$$(b * \beta(a)) * a' = (b * a) *a',$$
here we apply $\widetilde{\rm Axiom 2}$ for the set $(A \cup B)$ and
the fact $\beta(a) * a' = a* a'.$

For any unary operation $\omega \in \Omega'_1$,

$$\omega(\beta(b
\cdot a)) = \beta(\omega(b \cdot a)) = \beta(\omega(b) \cdot
\omega(a)),$$ here we apply condition 10 of Proposition 2.6,

$$\omega(b + \beta(a) - b ) = \omega(b) + \beta(\omega(a)) -
\omega(b).$$  As we have proved above these elements are equal.

Below we apply condition 11 of Proposition 2.6 and obtain:
$$\omega(\beta(b*a)) = \beta(\omega(b*a)) = \beta(\omega(b)*a),$$
$$\omega(b * \beta(a)) = \omega(b) * \omega(a).$$
As we have shown above these elements are equal. For $\omega_1,
\cdots, \omega_n$ the corresponding equalities are obtained
similarly.
\end{proof}

By Proposition 3.8, Definition 2.7 is equivalent to the following
one.

\begin{defi}
For any object $A$ in $\bC$ an actor of $A$ is an object
$\rm{Actor}(A)$, which acts on $A$ in $\bC$, and for any object
$C$ of $\bC$ and an action of $C$ on $A$ there is a unique
morphism $\vf:C\lra \rm{Actor}(A)$ with $c\cdot a=\vf(c)\cdot a$,
$c*a=\vf(c)*a$ for any $*\in\Om_2{}'$, $a\in A$ and $c\in C$.
\end{defi}

It is a well-known fact that the category of crossed modules in
the category of groups $\mathbb{X}Mod(\mathbb{G}r)$ is equivalent
to the category $\mathbb{G}$ with objects groups with the
additional two unary operations $\om_0, \om_1 : G \to G, G \in
\mathbb{G}r$ which are group homomorphisms satisfying conditions
\begin{enumerate}
\item[(1)] $\om_0 \om_1 = \om_1, \om_1 \om_0 = \om_0 $
 \item[(2)]
$\om_1(x)+y-\om_1(x)=  x +y -x,\ \  x,y \in {\rm Ker} \ \om_0$
\end{enumerate}
This category is a category of interest. The computations and
properties of actions in this category  and the direct checking of
identities (1), (2) show that $\cB(A)$ is an actor of $A \in
\mathbb{G}$. Thus the same is true for the category of crossed
modules $\mathbb{X}Mod(\mathbb{G}r)$. From the results of Norrie
\cite{No} it follows that the constructed by her the object
$A(T,G,\partial)$ for any crossed module $(T,G,\partial)$ is an
actor in the sense of Definition 2.7. Thus it follows that in the
category of interest $\mathbb{G}$ there exists an actor for any $A
\in \mathbb{G}$. By the Proposition 3.5 it follows that $\cB(A)$
is an actor for any $A \in \mathbb{G}$. It is another way of
proving that $\cB(A) =$ Actor$(A)$ in $\mathbb{G}$.

The category of precrossed modules is equivalent to the category
of interest $\mathbb{\bar G}$, which objects are groups with
additional two unary operations $\om_0, \om_1$, which are group
homomorphisms satisfying identity (1). By Theorem 3.7, $\cB(A) =$
GActor$(A)$, for any $A \in
 \mathbb{\bar G}$. It is easy to check that $\cB(A)$ satisfies
identity (1) and thus $\cB(A) \in \bar{\mathbb{G}}$, therefore
$\cB(A) =$ Actor$(A)$. From this we conclude that in the category
of precrossed modules always exists an actor.

As we have mentioned in the introduction internal object actions
were studied recently by F. Borceux, G. Janelidze and G. M. Kelly
\cite{BJK}, where the authors introduce a new notion of
representable action. From Theorem 6.3 of \cite{BJK}, applying
Proposition 3.8, it follows that in the case of category of
interest $\mathbb{C}$ the existence of representable object
actions is equivalent to the existence of an Actor$(A)$ for any $A
\in \mathbb{C}$ in the sense of Definition 2.7. Thus by Theorem
3.6, $\mathbb{C}$ has representable object actions if and only if
$\cB(A) \ltimes A \in \mathbb{C}$, for any $A \in \mathbb{C}$, and
if it is the case, the corresponding representing objects are
$\cB(A), A \in \mathbb{C}$.

\section{The case ${\bf \Om_2=\{+,*,*^\circ\}}$}

It is interesting to know in which kind of categories of interest
$\bC$ there exists  Actor$(A)$ for any object $A\in\bC$; or what
the sufficient conditions for the existence  of Actor$(A)$ for a
certain $A\in\bC$ are. In the case of groups ($\Om_2=\{+\}$), the
direct checking shows that $\cB(A)\in$ $\bG$r, and the action of
$\cB(A)$ on $A$ is a derived action. It follows also from
Propositions 3.1 and 3.2; thus $\cB(A)$ is an actor of $A$ by the
Theorem 3.6. This fact is also a consequence of Proposition 3.5
since it is well-known that Aut$(A)$ is an actor of $A$ in $\bG$r,
thus $\cB(A)\approx\text{Aut}(A)$. In the case of Lie algebras
($\Om_2=\{+,[\,,\,]\}$), the object $\cB(A)\in \bL$ie and the
action of $\cB(A)$ on $A$ is a derived action, so $\cB(A)$ is an
actor again in $\bL$ie and therefore $\cB(A)\approx\text{Der}(A)$.

Consider the case of Leibniz algebras. In this case we can define
the bracket operation for the elements of $\bB$ in two ways (see
Section 2 for the definition of the set $\bB$).

\begin{defi}
\begin{align*}
&\big[a,[\bbb_i,\bbb_j]\big]=\big[[a,b_i],b_j\big]-\big[[a,b_j],b_i\big],\\
&\big[[\bbb_i,\bbb_j],a\big]=\big[b_i,[b_j,a]\big]+\big[[b_i,a],b_j\big].
\end{align*}
\end{defi}

\begin{defi}
\begin{align*}
&\big[a,[\bbb_i,\bbb_j]\big]=\big[[a,b_i],b_j\big]-\big[[a,b_j],b_i\big],\\
&\big[[\bbb_i,\bbb_j],a\big]=-\big[b_i,[a,b_j]\big]+\big[[b_i,a],b_j\big].
\end{align*}
\end{defi}

The bracket operation $[b,b']$ for any $b,b'$ which are the
results of bracket  operations itself is defined according to
above formulas.

The addition is defined by:
\begin{align*}
&[\bbb_i+\bbb_j,a]=[b_i,a]+[b_j,a],\\
&[a,\bbb_i+\bbb_j]=[a,b_i]+[a,b_j].
\end{align*}
For any $b,b'\in \cB(A)$, $b+b'$ is defined by the same formulas.

The action of $\cB(A)$ on $A$ is defined according to Definition
4.1 or 4.2 respectively. So we have two different ways of
definition of an action. It is easy to check that non of them is
the derived action in $\bL$eib.

The algebras $\cB(A)$ defined by Definitions 4.1 and 4.2 are not
isomorphic.

\begin{cond}
For $A\in\bL$eib, and any two objects $B,C\in\bL$eib, which act on
$A$, we have
$$  \big[c,[a,b]\big]=-\big[c,[b,a]\big],   $$
$a\in A$, $b\in B$, $c\in C$.
\end{cond}

Note that in this condition under action we mean the derived
action. \vskip+2mm

\begin{example}
 If Ann$(A)=(0)$ or $[A,A]=A$, then $A$
satisfies Condition 1.
\end{example}
\begin{Prop}
For any object $A\in\bL${\rm eib}, the Definitions $4.1$ and $4.2$
give the same algebras if $A$ satisfies Condition $1$.
\end{Prop}

The proof follows directly from the definitions of operations in
$\cB(A)$ and Condition 1.

Below we mean that $\cB(A)$ is defined in one of the ways.

\begin{Prop}
For any $A\in\bL${\rm eib}, $\cB(A)$ is a Leibniz algebra. The set
of actions of $\cB(A)$ on $A$ is a set of derived actions if and
only if $A$ satisfies Condition $1$.
\end{Prop}

\begin{proof}
The computation shows that if Condition 1  holds then the same kind
of condition is fulfilled for $b, c \in \cB(A)$, from which follows
the result.
\end{proof}

\begin{theo}
For a Leibniz algebra $A$ there exists an actor if and only if $A$
satisfies Condition $1$. If it is the case, then
$\cB(A)=\rm{Actor}(A)$.
\end{theo}

\begin{proof}
By Proposition 4.3, $\cB(A)$ is always a Leibniz algebra and by
Theorem 3.7, $\cB(A)=G\text{Actor}(A)$. If $A$ satisfies Condition
1, by Proposition 4.4, $\cB(A)$ has a derived action on $A$ and thus
$\cB(A)=\text{Actor}(A)$. Conversely, if $A$ has an actor then
$\cB(A)=\text{Actor}(A)$ by Proposition 3.5, and so the action of
$\cB(A)$ on $A$ is a derived action, thus we have for any $a\in A$,
$b_i\in B_i$, $b_j\in B_j$, $i,j\in\bJ$ the following equalities
\begin{align*}
&\big[b_i,[a,b_j]\big]=\big[[b_i,a],b_j\big]-\big[[b_i,b_j],a\big],\\
&\big[b_i,[b_j,a]\big]=\big[[b_i,b_j],a\big]-\big[[b_i,a],b_j\big],
\end{align*}
from which follows Condition 1, which proves the theorem.
\end{proof}

We have an analogous picture for associative algebras. The
operations for the elements of $\bB$ (see Section 2 for the
notation) in this category are given by
\begin{equation}
\left.\begin{array}{l}
(\bbb_i*\bbb_j)*(a)=b_i*(b_j*a),\\
*(\bbb_i*\bbb_j)(a)=(a*b_i)*b_j,\\
(\bbb_i+\bbb_j)*(a)=b_i*a+b_j*a,\\
*(\bbb_i+\bbb_j)(a)=a*b_i+a*b_j.\\
\end{array}\right.
\end{equation}
The set of actions of $\cB(A)$ on $A$ is defined according to
(4.1).

\begin{cond}
For $A\in \mathbb{A}${\rm ss} and any two objects $B$ and $C$ from
$ \mathbb{A}${\rm ss} which have derived actions on $A$ we have
$$c*(a*b)=(c*a)*b,$$
for any $a\in A$, $b\in B$, $c\in C$.
\end{cond}

\begin{example} If Ann$(A)=(0)$ or $A^2=A$ then $A$
satisfies Condition 2. For this kind of associative algebras it is
proved in \cite{LL}  that $A\lra\text{Bim}(A)$ is a terminal
object in the category of crossed modules under~$A$.
 \end{example}
\begin{Prop}
For $A\in \mathbb{A}${\rm ss}, the algebra $\cB(A)$ is an
associative algebra and the set of actions of $\cB(A)$ on $A$
defined according to $(4.1)$ is the set of derived actions in $
\mathbb{A}${\rm ss} if and only if $A$ satisfies Condition ~$2$.
If it is the case $\cB(A)= \ ${\rm Actor}$(A)$.
\end{Prop}

The proof contains the analogous arguments as for the case of
Leibniz algebras and is left to the reader.

It is easy to see that in $ \mathbb{A}$ss and $\bL$eib generally
we have injections
$$  \cB(A)\lra\text{Bim}(A)\quad \text{and}\quad \cB(A)\lra\text{Bider}(A)\  $$
which are homomorphisms in $ \mathbb{A}$ss and $\bL$eib
respectively.

\begin{Prop}
Let $A$ be an associative algebra with the condition {\rm
Ann}$(A)=0$ or $A^2=A$. Then $\cB(A)\approx$ \ {\rm Bim}$(A)= \
${\rm Actor}$(A)$.
\end{Prop}

\begin{proof}
It is well-known that Bim$(A)$ is an associative algebra \cite{Mc}.
The action of Bim$(A)$ on $A$ (see Section 2) is not a derived
action in general, and the condition which fails is
\begin{equation}
f*(a*f')=(f*a)*f'
\end{equation}
for any $f=(f*,*f)$ and $f'=(f'*,*f')$ from Bim$(A)$. The direct
checking shows that in case Ann$(A)=(0)$ or $A^2=A$ identity (4.2)
holds for the action  \cite{LL}. For  any action of the object $B$
on $A$, $B\in \mathbb{A}$ss, we define $\vf:B\lra\text{Bim}(A)$ by
$\vf(b)=(b*,*b)$, which is a unique homomorphism with the property
that $\vf(b)*a=b*a$, $*\in\Om_2'$, since in Bim$(A)$ for any two
elements $f,f'\in\text{Bim}(A)$ from $f=f'$ follows that $f*=f'*$,
$*f=*f'$. Thus Bim$(A)$ is an actor of $A$ in $ \mathbb{A}$ss and
the isomorphism $\cB(A)\approx\text{Bim}(A)$ follows from
Proposition 3.5.
\end{proof}

We have the analogous result for Leibniz algebras.

\begin{Prop}
Let $A\in\bL${\rm eib} and {\rm Ann}$(A)=(0)$ or $[A,A]=A$. Then
$\cB(A)\approx$ \ {\rm Bider}$(A)=$ \ {\rm Actor}$(A)$.
\end{Prop}

\begin{proof}
We will follow the first definition of the bracket operation in
Bider$(A)$ (see Section 2, (2.5.1), (2.5.2)). The direct checking
shows that Bider$(A)$ is a Leibniz algebra (see Remark below and cf.
\cite{Lo}). The action of Bider$(A)$ on $A$ is not a derived action,
fails the following condition
\begin{equation}
\big[\vf,[a,\vf']\big]=\big[[\vf,a],\vf'\big]-\big[[\vf,\vf'],a\big],
\end{equation}
where $\vf=[[\;,\vf],[\vf,\;]]$ and $\vf'=[[\;,\vf'],[\vf', \
]]\in\text{Bider}(A)$.

From (2.5.2) we have
\begin{equation}
\big[\vf,[\vf',a]\big]=\big[[\vf,\vf'],a\big]-\big[[\vf,a],\vf'\big].
\end{equation}
We shall show that if Ann$(A)=(0)$ then
$[\vf,[\vf',a]]=-[\vf,[a,\vf']]$, and from (4.4) will follow (4.3).

Note that under the annulator we mean both left and right annualtor.
For any $a'\in A$ we have the following equalities:
\begin{align*}
\big[a',[\vf,[\vf',a]]\big]&=-\big[a',[[\vf',a],\vf]\big]=-\big[[a',[\vf',a]],\vf\big]
+\big[[a',\vf],[\vf',a]\big],\\
\big[a',[\vf,[a,\vf']]\big]&=-\big[[a',[a,\vf']],\vf\big]+\big[[a',\vf],[a,\vf']\big]\\
&=\big[[a',[\vf',a]],\vf\big]-\big[[a',\vf],[\vf',a]\big].
\end{align*}
Thus we obtain that for $a'\in A$
$$  \big[a',[\vf,[\vf',a]]+[\vf,[a,\vf']]\big]=0   \, .$$
In analogous way we show that
$$  \big[[\vf,[\vf',a]]+[\vf,[a,\vf']],a'\big]=0    \, .$$
From which we conclude that
$$  \big[\vf,[\vf',a]\big]+\big[\vf,[a,\vf']\big]=0 \, . $$
The case $[A,A]=A$ can be proved analogously. Thus we have a derived
action of Bider$(A)$ on $A$ and the crossed module
$A\lra\text{Bider}(A)$ ($a\mapsto([\;,a],[a,\;]$) has the universal
property of the actor object. By Proposition 3.5
$\cB(A)\approx\text{Bider}(A)$ which ends  the proof.
\end{proof}

\begin{remark}
 As we have also mentioned in Section 2, if
$[\vf,[\vf',a]]=-[\vf,[a,\vf']]$ for any $\vf=([\;,\vf],[\vf,\;])$
and $\vf'=([\;,\vf'],[\vf',\ ])$ from Bider$(A)$, then two
definitions of Bider$(A)$ according to (2.5.1), (2.5.2) and
(2.5.1), (2.5.2') coincide and this algebra is isomorphic to the
Leibniz algebra of biderivations defined by Loday \cite{Lo}.
\end{remark} \vskip+2mm

 In the category of $R$-modules over some ring
$R$, it is obvious that Actor($A$) = 0 for any $A$ since every
action is trivial in this category. The same result gives our
construction, $\cB(A)=0$ for any $R$-module $A$.

As it is in the case of associative algebras, in the category of
commutative associative algebras the condition for the action $(b_1
a) b_2 = b_1 (a b_2)$ fails; also in this category we must have $b a
= a b$, for $b \in \cB(A)$, and $b_1 b_2 = b_2 b_1$ for $b_1, b_2
\in \cB(A)$. All these conditions are satisfied and we have $\cB(A)
= $Actor($A$) in commutative associative algebras if and only if $A$
satisfies Condition 2. If Ann$(A)=(0)$  or $A^2=A$ then $A$
satisfies Condition 2. For this kind of commutative algebras,
Actor($A$) = Bim($A$) = M($A$) where M($A$) is the set of
multiplications (or multipliers) of $A$ \cite{LS,LL}, i.e.,
$k$-linear maps  $f:A \to A$ with $f(aa')=f(a)a'$.

In the category of alternative algebras Actor($A$) does not exist
for any $A$. The condition of Proposition 4.9 given below is not
fulfilled. If $A$ satisfies the condition $b(a c) = (b a) c$, for
$b, c \in B, C,$ respectively, where $B, C$ are any actions on $A$,
then $\cB(A)$ is an actor of $A$. Note that the condition given
above implies that $A$ is an associative algebra. So we have actors
for associative algebras with this condition (i.e. Condition 2 in
alternative algebras) in the category $\mathbb{A}$lt. We can
consider the weaker condition $(xy)z- x( yz)=z(y x) -(z y)x, x,y,z
\in A \cup (\cup B_i)$, which does not imply in general
associativity of $A$. This condition is important for  the
fulfilment of the action conditions for the elements from $\cup
B_i$. The existence of an actor in $\mathbb{A}$lt under this
condition can be studied in the future.
\bigskip

 Consider now a more general case, where $(\bC,\bE,\Om)$ is
a category of interest and $\Om_2=\{+,*,*^\circ\}$. In this case
Axiom 2 contains two identities
\begin{align*}
(y*z)*x&=W_1(y,z;x;*,*),\\
(y*z)*^\circ x&=W_2(y,z;x;*,*^\circ).
\end{align*}
So under Axiom 2 we mean the two identities above. $\widetilde{\rm
Axiom} 2$ will denote the corresponding identities from $\wt{\bE}$
(see Section 2). For $C\in\bC$ and $x,y,z\in C$ denote
$$  T=\big\{(y*(x*z),z*(x*y)),(y*(z*x),z*(y*x)),$$ $$((y*x)*z,(z*x)*y),((x*y)*z,(x*z)*y)\big\}.  $$

\begin{cond}

a) The words $W_1$ and $W_2$ in Axiom 2 contain at least one
element from each pair of the set $T$ so that each one can be
expressed by $W_1$ or $W_2$ (e.g. computing from Axiom 2 we must
have $y*(x*z) = W_2(x,z;y;*,*^{\circ})$, the analogous equalities
for other elements from $T$) in a direct way, i.e. not due to
identities from $\mathbb{E}$ or their consequences.

b) $W_1(x,y;z*t;*,*)$ and $W_2(z,t; x*y;*, *^{\circ})$ are the same
words up to commutativity of ``juxtapositions''.
\end{cond}

Here under ``juxtaposition'' we mean e.g. $x_1(x_2x_3)$ and thus
each member from the eight members in $W(\, )$.

\noindent  $\widetilde{\textbf{ Condition \ 3.}}$  \medspace
 It is
 analogous to Condition 3 but is stated for elements of certain
 $A\in  \mathbb{C}$ and the elements of its different actions,
 thus for $x,y,z,t \in A \cup(\cup B_i)$, whenever they have a
 sense. We admit that a) and b) conditions are fulfilled not
 necessarily in a direct way. So there can be applied identities
 from $\widetilde{\mathbb{E}}$ and the special properties of $A$ itself.

\begin{cond}
The final decompositions of the words $W(x*y, z;t;*,\overline{*})$
and $W(x,y;z;*,*) \overline{*}t, W(x*y,z;t;*^{\circ},\overline{*})$
and $W(x,y;z;*,*^{\circ}) \overline{*} \ t$ are the same  up to
commutativity of ``juxtapositions'', $\overline{*} \in
\{*,*^{\circ}\}$. We mean the corresponding indices for $W$ in each
case.
\end{cond}

\noindent  $\widetilde{\textbf{ Condition \ 4.}}$  \medspace It is
analogous to Condition 4 but we have $x,y,z,t \in A \cup(\cup B_i)$
and we mean that we have the equalities between pairs of words given
in Condition 4, applying identities from $\mathbb{E}$ and the
special properties of $A$.

\begin{Prop}
i) If $\cB(A)$ has the derived action on $A$ then $A$ satisfies
$\widetilde{Condition \ 3}$.

ii) If $\cB(A) \in \mathbb{C}$, then $A$ satisfies
$\widetilde{Condition \ 4}$.
\end{Prop}

\begin{proof}
i) From the Theorem 2.5 and the definition of the algebra $B \ltimes
A$ it follows that the conditions for derived actions  for $*$
operation follow from the Axioms 1 and 2, where $x_1, x_2, x_3 \in A
\cup B$, whenever it has a sense. The result now follows directly
from the conditions of the Proposition.

ii) It is obvious.
\end{proof}

\begin{remark}  Leibniz algebras, associative algebras,
satisfying Conditions 1, 2, respectively, are examples of i) in
Proposition 4.9. Note that identities in $\mathbb{E}$ involving only
once the operation $*$, e.g. $x*y = y*x$ or $x*y = - y*x$ play an
important role. This is the case e.g. of commutative associative and
Lie algebras.
\end{remark}
\bigskip

These conditions are not generally sufficient since we do not know
what kind of identities we have in $\mathbb{E}$. These conditions
usually can be not sufficient even in the case where
$\bE=\bE_G\cup\{{\rm Axiom 1,Axiom 2}\}$, since it may happen that
$\widetilde{\rm Condition \ 3}$ is not fulfilled when certain $b_i
\in B_i$ is replaced by the element of $b_i*b_p \in \cB(A)$ in the
identities involved in $\widetilde{\rm Condition \ 3}$. The same
note we can make concerning $\widetilde{\rm Condition \ 4}$.
\bigskip

Below we summarize for the case $\bE=\bE_G\cup\{{\rm Axiom 1,Axiom
2}\}$ our results and obtain

\begin{theo}
Let $\bC$ be a category of interest, where $\Om_2=\{+,*,*^\circ\}$
and $\bE=\bE_G\cup\{{\rm Axiom 1,Axiom 2}\}$.

{\rm a)} If $\mathbb{C}$ satisfies  Condition 3, then $\cB(A)$ has
a derived action on $A$ for any $A \in \mathbb{C}$.

{\rm b)} If $\mathbb{C}$ satisfies Condition 4, then $\cB(A)\in
\mathbb{C}$ for any $A \in \mathbb{C}$.

{\rm c)} If $\mathbb{C}$ satisfies Conditions 3 and 4, then
$\cB(A) = {\rm Actor}(A)$ for any $A \in \mathbb{C}$.

\end{theo}

If $\rm{\widetilde{Conditions \ 3}}$ and $4$ are satisfied for any
$A \in \mathbb{C}$, then these conditions are satisfied for all free
algebras;  it can involve certain identities which are consequences
of  $\widetilde{\rm Axiom 2}$. But these identities can be not true
for the elements of $A \cup \cB(A)$. In the case, where we do not
have  consequence identities of Axiom 2 from the fulfilment of
$\rm{\widetilde{Conditions \ 3}}$ and $4$  for any $A$ it follows
that Conditions 3 and 4 are also satisfied. Thus in this case
$\rm{\widetilde{Conditions \ 3}}$ and $4$ are sufficient conditions
for the existence of an actor. But it is important to note that if
Axiom 2 has no consequence identities, then $\widetilde{\rm
Conditions}$ 3 and 4 are always satisfied too. Actually we obtain
simpler conditions for this special case.

\begin{theo}
Let $\Om_2=\{+,*,*^\circ\}$, $\bE=\bE_G\cup\{{\rm Axiom 1,Axiom
2}\}$, and Axiom 2 does not imply new identities. $\cB(A) = {\rm
Actor}(A)$ for any $A \in \mathbb{C}$ if and only if $W_1, W_2$
contain at least one element from each pair of the set $T$.

\end{theo}

Below we consider the algebras with additional commutativity $(x *
y = y * x)$ or anticommutativity $(x * y = - y * x)$ condition on
the binary operation $*$. We will write $\mathbb{E} = \mathbb{E}_G
\cup \{{\rm Axiom \ 1}, {\rm Axiom \ 2}, {\rm (A)Comm}\}$. In the
corresponding category of interest $\mathbb{C}$, our construction
$\cB(A)$ must satisfy also (a)commutativity condition. For this
category we apply weaker forms of Conditions 3 and 4. We require
that they are fulfilled in a direct way using only
(a)commutativity property of the $*$ operation. In this case
(a)commutativity of $*$ operation in $\cB(A)$ guarantees the
identity

$$W_2(y,z;x;*,*^{\circ}) = (-)W_2(z,y;x;*,*^{\circ})
\eqno(4.5)$$

\noindent which must be fulfilled applying only (a)commutativity of
the $*$ operation and  commutativity of ``juxtapositions''. Note
that e.g. for commutative associative algebras (4.5) does not hold
in the way it is required above. For the corresponding equality in
this case we apply not only commutativity of the multiplication but
also associativity, thus Axiom 2 for this case.

\begin{theo}
Let $\mathbb{C}$ be a category of interest, $\Om_2=\{+,*\},
\mathbb{E}=\mathbb{E}_G \cup \{{\rm Axiom \ 1}, {\rm Axiom \ 2} ,
{\rm (A)Com}\}$. If {\rm Axiom  2}  does not imply new identities
and (4.5) holds, then $\cB(A)$ is an actor of A for any $A \in
\mathbb{C}$ if and only if $W( \ )$ in {\rm Axiom 2} contains at
least one element from each pair of the set T.
\end{theo}

\begin{proof}
Direct checking of identities.
\end{proof}

\begin{example}If Axiom 2 has the form $$x*(y*z) = -
y*(z*x) - z*(x*y) \eqno(4.6)$$ then the category $\mathbb{C}$ with
$\Om_2=\{+,*\}, \mathbb{E}=\mathbb{E}_G \cup \{{\rm Axiom \ 1},
{\rm Axiom \ 2} , {\rm Com}\}$ satisfies the conditions of Theorem
4.12. The same is true for the category $\mathbb{C}$ with the same
Axiom 2 and anticommutativity  property.
\end{example}

Note that (4.6) is equivalent to Jacobi identity, but the addition
is not commutative in our case.

 \

\end{document}